\documentclass[12pt]{article}%%%%%%%%%%%%%%%%%%%%%%%%%%%%%%%%%%%%%%%%%%%%%%%%%%%%%%%%%%%%%%%%%%%%%%%%%%%%%%%%%%%%%%%%%%%%%%%%%%%%%%%%%%%%%%%%%%%%%%%%%%%%%%%%%%%%%%%%%%%%%%%%%%%%%%%%%%%%%%%%%%%%%%%%%%%%%%%%%%%%%%%%%%%%%%%%%%%%%%%%%%%%%%%%%
\usepackage{amsfonts}
\usepackage{subfigure}
\usepackage{latexsym,amsmath,amssymb}
\usepackage{graphics,epstopdf}
\usepackage[pdftex]{graphicx}

\newtheorem{definition}{Definition}[section]
\newtheorem{theorem}[definition]{Theorem}

\newtheorem{remark}[definition]{Remark}

\numberwithin{equation}{section}

\numberwithin{equation}{section}

\begin{document}

\bigskip

\bigskip

\begin{center}

{\Large \textbf{Statistical approximation by $(p,q)$-analogue of Bernstein-Stancu Operators}}

\bigskip

\textbf{Asif Khan} and
 \textbf{Vinita Sharma}

Department of\ Mathematics, Aligarh Muslim University, Aligarh--202002, India%
\\[0pt]

asifjnu07@gmail.com; vinita.sha23@gmail.com \\[0pt%
]

\bigskip

\bigskip

\textbf{Abstract}
\end{center}
{In this paper, some approximation properties of $(p,q)$-analogue of Bernstein-Stancu Operators has been studied. Rate of statistical convergence by means of modulus of continuity and Lipschitz type maximal functions has been investigated. Monotonicity of $(p,q)$-Bernstein-Stancu Operators  and a global approximation theorem by means of Ditzian-Totik modulus of smoothness is established. A quantitative Voronovskaja type theorem is developed for these operators.
	Furthermore, we show comparisons and some illustrative graphics for the convergence of operators to a function}.\\



\section{Introduction and preliminaries}
\hspace{8mm} Mursaleen et al. \cite{mka1}  first applied the concept of $(p,q)$-calculus in approximation theory and introduced the $(p,q)$-analogue of Bernstein operators. Later on, based on $(p,q)$-integers, some approximation results for  Bernstein-Stancu operators, Bernstein-Kantorovich operators, Bleimann-Butzer and Hahn operators, $(p,q)$-Lorentz operators, Bernstein-Shurer operators, $(p,q)$-analogue of divided difference and Bernstein operators  etc.  have also been introduced by them in \cite{mur8,mka3,mka5,m4,zmn, mnfa}.\\

For similar works in approximation theory \cite{pp} based on $q$ and $(p,q)$-integers, one can refer \cite{acar1,acar2,acar3,acar4,aral, cai,ali,lupas,ma1,sofia,kang,kadak2,kac,mahmudov1,wafi}.

Motivated by the work of Mursaleen et al \cite{mka1}, the idea of $(p,q)$-calculus and its importance.\\

Very recently, Khalid et al. \cite{khalid1,khalid2,khalid3,kblossom} has given a nice application in computer-aided geometric design and applied these Bernstein basis  for construction of $(p,q)$-B$\acute{e}$zier curves and surfaces based on  $(p,q)$-integers which is further generalization of $q$-B$\acute{e}$zier curves and surfaces \cite{bezier,hp,pl,phillips,pp,lp}.  For similar works, one can refer \cite{bezier,hp}. Another advantage of using the parameter $p$ has been discussed in \cite{m4}.\\

Let us recall certain notations of $(p,q)$-calculus .\\

For any $p>0$ and $q>0,$ the $(p,q)$ integers $[n]_{p,q}$ are defined by

\begin{equation*}
[n]_{p,q}=p^{n-1}+p^{n-2}q+p^{n-3}q^2+...+pq^{n-2}+q^{n-1}\\
=\left\{
\begin{array}{lll}
\frac{p^{n}-q^{n}}{p-q},~~~~~~~~~~~~~~~~\mbox{when $~~p\neq q \neq 1$  } & \\
&  \\
n~p^{n-1},~~~~~~~~~~~~~~\mbox{ when $p=q\neq1$  } & \\
&  \\

[n]_q ,~~~~~~~~~~~~~~~~~~~\mbox{when $p=1$  }& \\

n ,~~~~~~~~~~~~~~~~~~~~~\mbox{ when $p=q=1$  }
\end{array}%
\right.
\end{equation*}

~where  $[n]_q $ denotes the $q$-integers and $n=0,1,2,\cdots$.\\

Obviously, it may be seen that $[n]_{p,q}= p^{n-1}[n]_{\frac{q}{p}}.\\$

The $(p,q)$-factorial is defined by
$$[0]_{p, q}!:=1~~\text{and}~~[n]!_{p, q}=[1]_{p, q}[2]_{p, q}\cdots [n]_{p, q}~~\text{if}~~n\ge 1.$$
Also the $(p,q)$-binomial coefficient is defined by
$${n \brack k}_{p, q}=\frac{[n]_{p, q}!}{[k]_{p, q}!~[n-k]_{p, q}!}~~\text{for all}~~n, k\in \mathbb N~~\text{with}~~n\ge k.$$

The formula for $(p,q)$-binomial expansion is as follows:
\begin{equation*}
(ax+by)_{p,q}^{n}:=\sum\limits_{k=0}^{n}p^{\frac{(n-k)(n-k-1)}{2}}q^{\frac{k(k-1)}{2}}
\left[
\begin{array}{c}
n \\
k%
\end{array}%
\right] _{p,q}a^{n-k}b^{k}x^{n-k}y^{k},
\end{equation*}

$$(x+y)_{p,q}^{n}=(x+y)(px+qy)(p^2x+q^2y)\cdots (p^{n-1}x+q^{n-1}y),$$
$$(1-x)_{p,q}^{n}=(1-x)(p-qx)(p^2-q^2x)\cdots (p^{n-1}-q^{n-1}x),$$\\

Details on $(p,q)$-calculus can be found in \cite{mah,jag,mka1,khalid1,khalid2,kblossom}.\\

The $(p,q)$-Bernstein Operators introduced by Mursaleen et al. for $0<q<p\leq 1$  in \cite{mka1} are as follow:
\begin{equation}\label{ee1}
B_{n,p,q}(f;x)=\frac1{p^{\frac{n(n-1)}2}}\sum\limits_{k=0}^{n}\left[
\begin{array}{c}
n \\
k%
\end{array}%
\right] _{p,q}p^{\frac{k(k-1)}2}x^{k}\prod\limits_{s=0}^{n-k-1}(p^{s}-q^{s}x)~~f\left( \frac{%
	[k]_{p,q}}{p^{k-n}[n]_{p,q}}\right) ,~~x\in \lbrack 0,1].
\end{equation}

Note when $p=1,$ $(p,q)$-Bernstein Operators given by \eqref{ee1} turns out to be $q$-Bernstein Operators.\\

Also, we have
\begin{align*}
(1-x)^{n}_{p,q}&=\prod\limits_{s=0}^{n-1}(p^s-q^{s}x) =(1-x)(p-qx)(p^{2}-q^{2}x)...(p^{n-1}-q^{n-1}x)\\
&=\sum\limits_{k=0}^{n} {(-1)}^{k}p^{\frac{(n-k)(n-k-1)}{2}} q^{\frac{k(k-1)}{2}}\left[
\begin{array}{c}
n \\
k% x^{k}
\end{array}%
\right] _{p,q}x^{k}
\end{align*}

Motivated by the above mentioned work on $(p,q)$-approximation and its application, this paper is organized as follows:
In Section 2, some basic results  for  $(p,q)$-analogue of Bernstein-Stancu Operators as given in \cite{mka3} has been recalled and based on it, second order moment is computed. In section 3, Korovkin's type statistical approximation properties has been studied for these operators. In section 4, rate of statistical convergence by means of modulus of continuity and Lipschitz type maximal functions has been investigated. Section 5 is based on monotonicity of $(p,q)$-Bernstein-Stancu Operators. In section 6, a global approximation theorem by means of Ditzian-Totik modulus of smoothness and a quantitative Voronovskaja type theorem is established.

The effects of the  parameters $p$ and $q$ for the convergence of operators to a function is shown in section 7 .\\

\section{$(p,q)$- Bernstein Stancu operators}

Mursaleen et. al in  \cite{mka3} introduced $(p,q)$-analogue of Bernstein-Stancu operators as follow:

\begin{equation}\label{ee2}
S_{n,p,q}(f;x)=\frac1{p^{\frac{n(n-1)}2}}\sum\limits_{k=0}^{n}\left[
\begin{array}{c}
n \\
k%
\end{array}%
\right] _{p,q}p^{\frac{k(k-1)}2}x^{k}\prod\limits_{s=0}^{n-k-1}(p^{s}-q^{s}x)~~f\left( \frac{%
	{p^{n-k}[k]_{p,q}+\alpha}}{[n]_{p,q}+\beta}\right) ,~~x\in \lbrack 0,1].
\end{equation}
where $\alpha$ and $\beta$ are real numbers which satisfy $0\leq\alpha \leq \beta$.\\

Note that for $\alpha=\beta=0,$ $(p,q)$-Bernstein-Stancu operators given by \eqref{ee2} reduces into $(p,q)$-Bernstein operators as given in \cite{mka1}.\\

Also for $p=1$, $(p,q)$-Bernstein-Stancu operators given by \eqref{ee2} turn out to
be $q$-Bernstein-Stancu operators.\\

For $p=q=1,$ it reduces to classical Bernstein-Stancu operators.\\

We have the following auxiliary lemmas:\newline

\parindent=0mm\textbf{Lemma 2.1.} For $x\in \lbrack [0,1],~0<q<p\leq 1$, and $\alpha,\beta \in\mathbb{R}$ with $0\leq\alpha \leq \beta$, we
have\\
\newline
(i)~~$S_{n,p,q}(1;x)=~1$,\newline
(ii)~$S_{n,p,q}(t;x)=~\frac{[n]_{p,q}}{[n]_{p,q}+\beta}x+\frac{\alpha}{[n]_{p,q}+\beta}$,\newline
(iii)~$S_{n,p,q}(t^{2};x)=\frac{q[n]_{p,q}[n-1]_{p,q}}{([n]_{p,q}+\beta)^2}x^2+\frac{[n]_{p,q}(2\alpha+p^{n-1})}{([n]_{p,q}+\beta)^2}x+\frac{\alpha^2}{([n]_{p,q}+\beta)^2}$.\newline

\textbf{Proof:} Proof is given in \cite{mka3}  using the identity
\begin{equation}
\sum\limits_{k=0}^{n}\left[
\begin{array}{c}
n \\
k%
\end{array}%
\right] _{p,q}p^{\frac{k(k-1)}2}x^{k}\prod\limits_{s=0}^{n-k-1}(p^{s}-q^{s}x)={p^{\frac{n(n-1)}2}}.
\end{equation}

We give complete  proof of Lemma 1 (iii)

(iii)
\begin{eqnarray*}
	S_{n,p,q}(t^2;x) &=&\frac1{p^{\frac{n(n-1)}2}}\sum\limits_{k=0}^{n}\left[
	\begin{array}{c}
		n \\
		k%
	\end{array}%
	\right] _{p,q}p^{\frac{k(k-1)}2}x^{k}\prod\limits_{s=0}^{n-k-1}(p^{s}-q^{s}x)~~{\bigg(\frac{p^{n-k}[k]_{p,q}+\alpha}{[n]_{p,q}+\beta}\bigg)}^2\\
	&=&\frac{1}{([n]_{p,q}+\beta)^2}~\frac1{p^{\frac{n(n-1)}2}}\Bigg[p^{2n}\sum\limits_{k=0}^{n}\left[
	\begin{array}{c}
		n \\
		k%
	\end{array}%
	\right] _{p,q}p^{\frac{k(k-1)}2}x^{k}\prod\limits_{s=0}^{n-k-1}(p^{s}-q^{s}x)~\frac{[k]_{p,q}^2}{p^{2k}}\\
	&&+2\alpha ~p^n\sum\limits_{k=0}^{n}\left[
	\begin{array}{c}
		n \\
		k%
	\end{array}%
	\right] _{p,q}p^{\frac{k(k-1)}2}x^{k}\prod\limits_{s=0}^{n-k-1}(p^{s}-q^{s}x)~\frac{[k]_{p,q}}{p^{k}}\\
	&&+\alpha^2~\sum\limits_{k=0}^{n}\left[
	\begin{array}{c}
		n \\
		k%
	\end{array}%
	\right] _{p,q}p^{\frac{k(k-1)}2}x^{k}\prod\limits_{s=0}^{n-k-1}(p^{s}-q^{s}x)\Bigg].\\
\end{eqnarray*}
$$S_{n,p,q}(t^2;x)=\frac{1}{([n]_p,q+\beta)^2} [(A)+(B)+(C)]$$

\begin{eqnarray*}
	(A)&=&\frac {1}{p^{\frac{n(n-1)}2}} p^{2n}\sum\limits_{k=0}^{n}\left[
	\begin{array}{c}
		n \\
		k%
	\end{array}%
	\right] _{p,q}p^{\frac{k(k-1)}2}x^{k}\prod\limits_{s=0}^{n-k-1}(p^{s}-q^{s}x)~\frac{[k]_{p,q}^2}{p^{2k}}\\
	&=&\frac{p^{2n}}{p^{\frac{n(n-1)}2}}\sum\limits_{k=0}^{n} \frac{[n]}{[k]}\left[
	\begin{array}{c}
		n-1 \\
		k-1%
	\end{array}%
	\right] x^{k}(1-x)^{n-k}~\frac{[k]^2}{p^{2k}}
\end{eqnarray*}
On shifting the limits and using $[k+1]_{p,q}=p^k+q[k]_{p,q}$, we get our desired result.
\begin{eqnarray*}
	(A)&=&\frac{p^{2n}}{p^{\frac{n(n-1)}2}}\sum\limits_{k=0}^{n-1} \left[
	\begin{array}{c}
		n-1 \\
		k%
	\end{array}%
	\right] x^{k}(1-x)^{n-k-1}~\frac{p^k+q[k]}{p^{2k+2}}\\
	&&= \frac{p^{2n-2}[n]x}{p^{\frac{n(n-1)}{2}}}\Bigg[ p^{\frac{(n-1)(n-2)}{2}}+ \frac{q[n-1]x}{p}\sum\limits_{k=0}^{n-2} \left[
	\begin{array}{c}
		n-2 \\
		k%
	\end{array}%
	\right] x^{k}(1-x)^{n-k-2}\Bigg]\\
	&&=p^n[n]x+q[n][n-1]x^2
\end{eqnarray*}

n\\Similarly
$$(B)=\frac{2\alpha ~p^n}{p^{\frac{n(n-1)}{2}}}\sum\limits_{k=0}^{n}\left[
\begin{array}{c}
n \\
k%
\end{array}%
\right] _{p,q}x^{k}(1-x)^{n-k}~\frac{[k]_{p,q}}{p^{k}}=2\alpha [n]x$$
and
$$(C)=\frac{\alpha^2}{p^{\frac{n(n-1)}{2}}}~\sum\limits_{k=0}^{n}\left[
\begin{array}{c}
n \\
k%
\end{array}%
\right] _{p,q}x^{k}(1-x)^{n-k}={\alpha}^2$$

\parindent=0mm\textbf{Lemma 2.2}. For $x\in [0,1],~0<q<p\leq 1$ and $\alpha,\beta \in\mathbb{R}$ with $0\leq\alpha \leq \beta$,\\

Let n be any given natural number, then
\begin{eqnarray*}
	S_{n,p,q}\bigl{(}(t-x)^2;x\bigl{)}&= \big\{ \frac{q[n]_{p,q}[n-1]_{p,q}-[n]_{p,q}^2+{\beta}^2}{([n]_{p,q}+\beta)^2}\big\}x^2+\big\{ \frac{p^{n-1}[n]_{p,q}-2 \alpha \beta}{([n]_{p,q}+\beta)^2} \big\}x+\frac{\alpha^2}{([n]_{p,q}+\beta)^2}\\
	&\leq\frac{[n]_{p,q}p^{n-1} - 2\alpha\beta}{2([n]_{p,q}+\beta)^2}\phi^2(x)\leq \frac{[n]_{p,q}}{[n]_{p,q}+\beta}\phi^2(x)\\
\end{eqnarray*}

\section{ Main Results}

\subsection { Korovkin type approximation theorem}

We know that $C[a,b]$ is a
Banach space with norm
\begin{equation*}
\Vert f\Vert _{C[a,b]}:=\sup\limits_{x\in \lbrack a,b]}|f(x)|,~f\in C[a,b].
\end{equation*}
For typographical convenience, we will write $\Vert .\Vert $ in place of $\Vert .\Vert _{C[a,b]}$ if no confusion arises.\newline

\begin{definition} {\em Let $C[a,b]$ be the linear space of all real valued continuous functions $f$
		on $[a,b]$ and let $T$ be a linear operator which maps $C[a,b]$ into itself.
		We say that $T$ is $positive$ if for every non-negative $f\in $ $C[a,b],$ we
		have $T(f,x)\geq 0$ for all $x\in $ $[a,b]$ .}
\end{definition}

\parindent=8mm The classical Korovkin type approximation theorem can be stated as follows \cite{brn,korovkin};\\

Let $T_n: C[a, b] \to C[a, b]$ be a sequence of positive linear operators.
Then $\lim_{n\to\infty}\|T_{n}(f; x)-f(x)\|_\infty=0,\,\,\textrm{for
	all}~f\in C[a, b]$ if and only if $\lim_{n\to\infty}\|T_{n}(f_{i}; x)-f_i(x)\|_\infty=0,\,\,\textrm{for each}\,\,i=0,1,2,$ where the test function $f_i(x)=x^i$.

In next section, we study a statistical approximation properties of the operator $S_{n,p, q}$.

\subsection{Statistical approximation}

The statistical version of Korovkin theorem for sequence of positive linear operators has been
given by Gadjiev and Orhan \cite{go39}.

%Also this type approximations have been studied by many authors the reader may refer to \cite{mur2,moh1,ohh}.

Let $K$ be a subset of the set $\mathbb{N}$ of natural numbers. Then, the asymptotic density $\delta(K)$ of $K$ is defined as $\delta(K)=\lim_{n}\frac{1}{n}\big|\{k\leq n~:~k \in K\}\big|$ and $|.|$ represents the cardinality of the enclosed set. A sequence $x=(x_k)$ said to be statistically convergent to the number $L$ if for each $\varepsilon >0$, the set $K(\varepsilon)=\{k\leq n:|x_k-L|>\varepsilon\}$ has asymptotic density zero (see \cite{erd37,fast}), i.e.,
\begin{eqnarray*}\label{117}
	\lim_{n}\frac{1}{n}\big|\{k\leq n:|x_k-L|\geq \varepsilon\}\big|=0.
\end{eqnarray*}
In this case, we write $st-\lim x =L$.

Let us recall the following theorem:\\

\begin{theorem}\label{ta}\cite{go39}
	Let $A_n$ be the sequence of linear positive operators from $C[0,1]$ to $C [0,1]$ satisfies the conditions
	
	$st-\lim\limits_{n}\|S_{n,p,q}(( t^\nu;x))- (x)^\nu \|_C[0,1] = 0 $ for $\nu = 0,~ 1,~ 2.$
	then for any function $f\in C[0,1],$
	$st-\lim\limits_{n} \|S_{ n,p,q}(f) - f \|_C[0,1] = 0.$
\end{theorem}

\subsection{Korovkin Type statistical approximation properties}

The main aim of this paper is to obtain the korovkin type statistical approximation properties of operators defined in (\ref{ee2}) with the help of Theorem (\ref{ta}).\\

\begin{remark}\label{r5.1}
	For $q\in(0,1)$ and $p\in(q,1]$, it is obvious that
	$\lim\limits_{n\to\infty}[n]_{p,q}=0 $ or $\frac1{p-q}$. In order to reach to convergence
	results of the operator $L^{n}_{p,q}(f;x),$ we take a sequence $q_n\in(0,1)$ and $p_n\in(q_n,1]$
	such that $\lim\limits_{n\to\infty}p_n=1,$ $\lim\limits_{n\to\infty}q_n=1$ and $\lim\limits_{n\to\infty}p_n^n=1,$ $\lim\limits_{n\to\infty}q_n^n=1$. So we get
	$\lim\limits_{n\to\infty}[n]_{p_n,q_n}=\infty$.
\end{remark}

\begin{theorem}
	Let $S_{n,p,q}$ be the sequence of operators and the sequence $ p=p_n$ and $q=q_n$ satisfying Remark $(\ref{r5.1})$  then for any function $f\in C[0,1]$\\
	$$ st-\lim\limits_{n} \|~S_{n,p_n,q_n}{(f,.)}-f\|=0$$
\end{theorem}

\textbf{Proof:}

Clearly for $\nu=0,$ $$ S_{n,p,q}{(1,x)}=1,$$

which implies $$ st-\lim\limits_{n}\|S_{n,p,q}(1;x)~-1~\|~~=~~0.$$\\

For $\nu~=~1$\\

\begin{align*}
\|S_{n,p,q}~(t;x)~-~x~~\|&\leq~\bigg|\frac{[n]_{p,q}}{[n]_{p,q}+\beta}x~~+~~\frac{\alpha}{[n]_{p,q}+\beta}~~-~~x\bigg|\\
&= \bigg|\bigg(\frac{[n]_{p,q}}{[n]_{p,q}+\beta}~~-~~1\bigg)x~+~\frac{\alpha}{[n]_{p,q}+\beta}\bigg|\\
&\leq \bigg|\frac{[n]_{p,q}}{[n]_{p,q}+\beta}~~-~~1\bigg|~~+~~\bigg|\frac{\alpha}{[n]_{p,q}+\beta}\bigg|.
\end{align*}

For a given $ \epsilon >0$, let us define the following sets.\\

$$U = \{n : \|S_{n,p,q}(t;x) -x\|\geq\epsilon\}$$\\
$$U^{\prime} = \{n: 1 - \frac{[n]_{p,q}}{[n]_{p,q}+\beta}\}  \geq \epsilon $$\\
$$U^{\prime\prime} = \{n: \frac{\alpha}{[n]_{p,q} +\beta}\geq\epsilon\}$$

%It is obvious that $U \subseteq U^{\prime\prime}\cup U^{\prime},$ and by using $(\ref{e6}),$ $$ st~~-~~\lim\limits_{n}\bigg(1 - \f0rac{[n]_{p,q}}{[n]_{p,q}+\beta}\bigg) =0.$$\\

So using  $\delta \{k\leq n:1-\frac{[n]_{p,q}}{[n]_{p,q}+\beta}~~\geq \epsilon\}, $

then we get
$$st-\lim\limits_{n}\|S_{n,p,q}(t;x) - x\|=0.$$

Lastly for $\nu=2,$ we have \\
\begin{align*}
\|S_{n,p,q}(t^2:x)- x^2\|&\leq \big|\frac{q[n]_{p,q}[n-1]_{p,q}}{{([n]_{p,q}+\beta)}^{2}}~~-1\big|\\
&+\big|\frac{[n]_{p,q}(2\alpha+p^{n-1})}{[n]_{p,q}+\beta}^{2}x\big|+\big|\frac{\alpha^2}{{([n]_{p,q}+\beta)}^{2}}\big|.
\end{align*}

If we choose

$$\alpha_n=\frac{q[n]_{p,q}[n-1]_{p,q}}{{([n]_{p,q}+\beta)}^{2}}~~-1$$\\
$$\beta_n=\frac{[n]_{p,q}(2\alpha+p^{n-1})}{[n]_{p,q}+\beta}^{2}$$\\
$$\gamma_n=\frac{\alpha^2}{{([n]_{p,q}+\beta)}^2}$$\\

$st-\lim\limits_{n}\alpha_n~~=~~st-\lim\limits_{n}\beta_n~~=~~st-\lim\limits_{n}\gamma_n~~=~~0$\\

Now given $\epsilon >0$, we define the following four sets:\\
$$U~~=~~\|S_{n,p,q}(t^2:x)- x^2\|\geq \epsilon$$\\

$$U_{1} =\{n:\alpha_{n} \geq \frac{\epsilon}{3}\}$$\\
$$U_{2}=\{n:\beta _{n}\geq \frac{\epsilon }{3}\}$$\\
$$U_{3}=\{n:\gamma _{n}\geq \frac{\epsilon }{3}\}.$$\\

It is obvious that$ U \subseteq U_1\bigcup U_2\bigcup U_3. $ Thus we obtain \\

$\delta\{K\leq n:\|S_{n,p,q}(t^2:x)- x^2\|\geq\epsilon\}$\\

$\leq\delta\{K\leq n:\alpha_{n} \geq \frac{\epsilon}{3} \}~+~\delta\{K\leq n:\beta _{n}\geq \frac{\epsilon }{3}\}+\delta\{K\leq n:\gamma _{n}\geq \frac{\epsilon }{3}\}$\\

So the right hand side of the inequalities is zero by $( \ref{117}).$\\

Then\\

$$st-\lim\limits_{n}\|S_{n,p,q}(t;x) -x\|=0$$ holds and thus  the proof is completed.

If we choose

$$\alpha_n=\frac{q[n]_{p,q}[n-1]_{p,q}}{{([n]_{p,q}+\beta)}^{2}}~~-1$$\\
$$\beta_n=\frac{[n]_{p,q}(2\alpha+p^{n-1})}{[n]_{p,q}+\beta}^{2}$$\\
$$\gamma_n=\frac{\alpha^2}{{([n]_{p,q}+\beta)}^2}$$\\

$st-\lim\limits_{n}\alpha_n~~=~~st-\lim\limits_{n}\beta_n~~=~~st-\lim\limits_{n}\gamma_n~~=~~0$\\

Now given $\epsilon >0$, we define the following four sets:\\
$$U~~=~~\|S_{n,p,q}(t^2:x)- x^2\|\geq \epsilon$$\\

$$U_{1} =\{n:\alpha_{n} \geq \frac{\epsilon}{3}\}$$\\
$$U_{2}=\{n:\beta _{n}\geq \frac{\epsilon }{3}\}$$\\
$$U_{3}=\{n:\gamma _{n}\geq \frac{\epsilon }{3}\}.$$\\

It is obvious that$ U \subseteq U_1\bigcup U_2\bigcup U_3. $ Thus we obtain \\

$\delta\{K\leq n:\|S_{n,p,q}(t^2:x)- x^2\|\geq\epsilon\}$\\

$\leq\delta\{K\leq n:\alpha_{n} \geq \frac{\epsilon}{3} \}~+~\delta\{K\leq n:\beta _{n}\geq \frac{\epsilon }{3}\}+\delta\{K\leq n:\gamma _{n}\geq \frac{\epsilon }{3}\}$\\

So the right hand side of the inequalities is zero by $( \ref{117}).$\\

Then\\

$$st-\lim\limits_{n}\|S_{n,p,q}(t;x) -x\|=0$$ holds and thus  the proof is completed.

\section{Rate of Statstical Convergence}

In this part, rates of statistical convergence of the operators $(\ref{ee2} )$ by means of modulus of continuity and LIPSCHITZ TYPE maximal functions are introduced.\\

The modulus of continuity for the space of function $ f\in C[0,1]$  is defined by\\

$$ w(f;\delta)=\sup\limits_{x,t\in C[0,1],~~ |t-x|<\delta} |f(t)-f(x)|$$\\

where ${w}(f;\delta)$ satisfies the following conditions:~~for all $f\in C[0,1],$\\

\begin{equation}\label{e118}
\lim\limits_{\delta\rightarrow0}~w (f;\delta) = 0.
\end{equation}
and
\begin{equation}\label{e119}
|f(t)-f(x)|\leq w(f;\delta)\bigg(\frac{|t-x|}{\delta}+ 1\bigg)
\end{equation}

\begin{theorem}
	
	Let the sequence $ p=p_n$ and $q=q_n$ satisfy for $0<q_n<p_n\leq1$, so we have \\
	
	$$|S_{n,p,q}(t;x) -f(x)|\leq w(f;\sqrt{\delta_n(x)})(1+q_n)$$\\
	
	where
	
	\begin{equation}\label{e10}
	\delta_n(x)=\frac{1}{([n]_{p,q}+\beta)^{2}}[(q[n]_{p,q}[n-1]_{p,q} -{[n]}^2+\beta^{2})x^2~~+~~([n]_{p,q}p^{(n-1)}-2\alpha\beta)x~+\alpha^2].
	\end{equation}
	
\end{theorem}

Proof: $|S_{n,p,q}(t;x) -f(x)|\leq S_{n,p,q}(|f(t)-f(x)|:x)$\\

by using $(\ref{e119}),$ we get \\

$$|S_{n,p,q}(t;x) -f(x)|\leq w(f;\delta)\{S_{n,p,q}(1;x)+\frac{1}{\delta}S_{n,p,q}(|t-x|:x)\}.$$\\

By using Cauchy Schwarz inequality, we have\\

\begin{align*}
|S_{n,p,q}(t;x)-f(x)|&\leq w(f;\delta_n)\bigg(1+\frac{1}{\delta_n}[(S_{n,p,q}(t-x)^2;x)]^{\frac {1}{2}}~~[S_{n,p,q}(1;x)]^{\frac {1}{2}}\bigg)\\
&\leq w(f;\delta_n)\bigg(1+\frac{1}{\delta_n}\bigg\{\frac{1}{([n]_{p,q}+\beta)^{2}}[(q[n]_{p,q}[n-1]_{p,q} -{[n]}^2\\~~
&+~~\beta^2)x^2~~+~~([n]_{p,q}p^{(n-1)}~-~2\alpha\beta)x~~
+\alpha^2]\bigg\}\bigg)
\end{align*}
so it is obvious by choosing $\delta_n$ as in $(\ref{e10})$ the theorem is proved.\\

Notice that by the condition in (\ref{e118}) $st-\lim\limits_{n} \delta_n =0,$ by $(\ref{e118})$ we have \\
$$st-\lim\limits_{n} w(f;\delta)~~=~~0.$$\

This gives us the pointwise rate of statistical convergence of the operators $S_{n,p,q}(f;x)~ \text{to}~ f(x).$\\

\section{Monotonicity for convex functions}
Oru\c{c} and Phillips proved that when the function $f$ is convex on $[0,1]$, its $q$-Bernstein operators are monotonic decreasing. In this section we will study the monotonicity of $(p,q)$-Bernstein Stancu operators.\\
\begin{theorem}
	If f is convex function on $[0,1],$ then
	$S_{n,p,q}(f;x)\geq f(x), $    $0\leq x \leq1$\\
	for all $n\geq 1$ and $0 < q < p \leq 1$
\end{theorem}

\textbf{Proof:} We consider the knots  $x_k = \frac{p^{n-k}[k]_{p,q}}{[n]_{p,q}},$

$$\lambda_k = \left[
\begin{array}{c}
n \\
k%
\end{array}%
\right] _{p,q} p^{\frac{k(k-1)-n(n-1)}{2}} x^k \prod\limits_{s=0}^{n-k-1}(p^{s}-q^{s}x),~~~~
0 \leq k \leq n.$$

Using Lemma 2.1, it follows that \\

$$\lambda_0+\lambda_1+\lambda_2+................\lambda_n  =  1$$\\
$$x_0\lambda_0+x_1\lambda_1+x_2\lambda_2+................x_n\lambda_n  =  x.$$\\

From the convexity of the function $f,$ we get\\

$S_{n,p,q}(f;x) =\sum\limits_{k=0}^{n}\lambda_k f(x_k)\geq f\bigg(\sum\limits_{k=0}^{n}\lambda_k x_k\bigg ) =f(x).$\\

\begin{theorem}
	Let f be convex on $[0,1]$. Then$S_{{n-1},p,q}(f;x)\geq S_{n,p,q}(f;x)$ for $0 < q < p\leq 1, $ $0 \leq x \leq 1,$ and $ n \geq 2 $. If $f \in C[0,1]$ the inequality holds strictly for $0 < x < 1 $ unless f is linear in each of the intervals between consecutive knots $\frac{p^{n-k-1}[k]_{p,q}}{[n]_{p,q}}$, $0 \leq k\leq n-1 $, in which case we have the equality.\\
\end{theorem}

\textbf{Proof:}
For $0<q<p\leq1,$ we begin by writing\\

$$\prod\limits_{s=0}^{n-1}(p^s -q^sx)^{-1} [S_{n-1,p,q}(f;x) - S_{n,p,q}(f;x)]$$\\

\begin{eqnarray*}
	&=&\prod\limits_{s=0}^{n-1}(p^s -q^sx)^{-1}\bigg [\sum\limits_{k=0}^{n-1} \left[
	\begin{array}{c}
		n-1 \\
		k%
	\end{array}%
	\right] _{p,q} p^{\frac{k(k-1)-(n-2)(n-1)}{2}} x^k \prod \limits_{s=0}^{n-k-2}(p^s -q^sx)f \bigg (\frac {p^{n-k-1}[k]_{p,q}+ \alpha}{[n]_{p,q}+ \beta}\bigg)\\
	&&-\sum\limits_{k=0}^{n} \left[
	\begin{array}{c}
		n \\
		k%
	\end{array}%
	\right] _{p,q} x^k p^{\frac{k(k-1)-n(n-1)}{2}} \prod\limits_{s=0}^{n-k-1}(p^s -q^sx)f\bigg(\frac {p^{n-k}[k]_{p,q}+ \alpha}{[n]_{p,q}+ \beta}\bigg)\bigg]\\
	&=&\sum\limits_{k=0}^{n-1} \left[
	\begin{array}{c}
		n-1 \\
		k%
	\end{array}%
	\right] _{p,q} p^{\frac{k(k-1)-(n-2)(n-1)}{2}} x^k \prod\limits_{s=n-k-2}^{n-1}(p^s -q^sx)^{-1}f\bigg (\frac {p^{n-k-1}[k]_{p,q}+ \alpha}{[n]_{p,q}+ \beta}\bigg)\\
	&&-\sum\limits_{k=0}^{n}\left[
	\begin{array}{c}
		n \\
		k%
	\end{array}%
	\right] _{p,q} x^k p^{\frac{k(k-1)-n(n-1)}{2}}\prod\limits_{s=n-k-1}^{n-1}(p^s -q^sx)^{-1} f\bigg (\frac {p^{n-k-1}[k]_{p,q}+ \alpha}{[n]_{p,q}+ \beta}\bigg).
\end{eqnarray*}

Denote\\
\begin{equation}\label{e14}
\psi_{k}(x) = p^{\frac{k(k-1)}{2}} x^k \prod\limits_{s=n-k-1}^{n-1}(p^s -q^sx)^{-1}
\end{equation}
and using the following relation:\\

\begin{align*}
p^{n-1} p^{\frac{k(k-1)}{2}} x^k \prod\limits_{s=n-k-1}^{n-1}(p^s -q^sx)^{-1} = p^k \psi_{k}(x)+q^{n-k-1}\psi_{k+1}(x).\\
\end{align*}

We find\\

\begin{equation*}
\prod\limits_{s=0}^{n-1}(p^s -q^sx)^{-1}[S_{n-1,p,q}(f;x) - S_{n,p,q}(f;x)]\\
\end{equation*}
\begin{eqnarray*}
	&=&\sum\limits_{k=0}^{n-1}\left[
	\begin{array}{c}
		n-1\\
		k%
	\end{array}%
	\right] _{p,q} p^{\frac{-(n-2)(n-1)}{2}} p^{-(n-1)}(p^k \psi_{k}(x)+q^{n-k-1}\psi_{k+1}(x)) f\bigg (\frac {p^{n-k-1}[k]_{p,q}+ \alpha}{[n]_{p,q}+ \beta}\bigg)\\
	&&-\sum\limits_{k=0}^{n}\left[
	\begin{array}{c}
		n \\
		k%
	\end{array}%
	\right] _{p,q} p^{\frac{-n(n-1)}{2}}\psi_{k}(x)f\bigg (\frac {p^{n-k-1}[k]_{p,q}+ \alpha}{[n]_{p,q}+ \beta}\bigg)\\
	%\end{align*}
	%\begin{align*}
	&=&p^{\frac{-n(n-1)}{2}}\bigg[\sum\limits_{k=0}^{n-1}\left[
	\begin{array}{c}
		n-1\\
		k%
	\end{array}%
	\right]_{p,q} p^k \psi_{k}(x)f\bigg(\frac{p^{n-k-1}[k]_{p,q}+\alpha}{[n]_{p,q}+ \beta}\bigg)\\
	&&+ \sum\limits_{k=1}^{n}\left[
	\begin{array}{c}
		n-1 \\
		k-1
	\end{array}%
	\right]_{p,q} q^{n-k} \psi_{k}(x) f \bigg(\frac {p^{n-k}[k]_{p,q}+ \alpha}{[n]_{p,q}+ \beta}\bigg)-\sum\limits_{k=0}^{n}\left[
	\begin{array}{c}
		n \\
		k%
	\end{array}%
	\right]_{p,q} \psi_{k}(x) f \bigg(\frac {p^{n-k}[k]_{p,q} + \alpha}{[n]_{p,q} + \beta} \bigg)\bigg]\\
	%\end{align*}
	%\begin{align*}
	&=&p^{\frac{-n(n-1)}{2}} \sum \limits_{k=1}^{n-1}\Bigg\{\left[
	\begin{array}{c}
		n-1\\
		k%
	\end{array}%
	\right]_{p,q} p^k f \bigg(\frac {p^{n-k-1}[k]_{p,q}+ \alpha}{[n]_{p,q}+ \beta}\bigg)\\
	&&+\left[
	\begin{array}{c}
		n-1\\
		k-1%
	\end{array}%
	\right] _{p,q}  q^{n-k} f\bigg(\frac{p^{n-k}[k]_{p,q}+ \alpha}{[n]_{p,q}+ \beta}\bigg) -\left[
	\begin{array}{c}
		n \\
		k%
	\end{array}%
	\right]_{p,q} f \bigg (\frac {p^{n-k}[k]_{p,q}+ \alpha}{[n]_{p,q}+ \beta}\bigg)\Bigg\} \psi_{k}(x)\\
	%\end{align*}
	%\begin{align*}
	&=&p^{\frac{-n(n-1)}{2}} \sum\limits_{k=1}^{n-1}\left[
	\begin{array}{c}
		n\\
		k%
	\end{array}
	\right]_{p,q}\Bigg\{\frac {[n-k]_{p,q}}{[n]_{p,q}}p^k f\bigg(\frac {p^{n-k-1}[k]_{p,q}+ \alpha}{[n]_{p,q}+ \beta}\bigg)\\
	&&+\frac{[k]_{p,q}}{[n]_{p,q}} q^{n-k} f\bigg(\frac{p^{n-k}[k]_{p,q}+ \alpha}{[n]_{p,q}+ \beta}\bigg)-f\bigg(\frac {p^{n-k}[k]_{p,q}+ \alpha}{[n]_{p,q}+ \beta}\bigg)\Bigg\}\psi_{k}(x)\\
	%\end{align*}
	%\begin{align*}
	&=&p^{\frac{-n(n-1)}{2}} \sum \limits_{k=1}^{n-1}\left[
	\begin{array}{c}
		n\\
		k%
	\end{array}
	\right]_{p,q} a_k \psi_{k}(x)
\end{eqnarray*}

where\\

$a_k =\frac{[n-k]_{p,q}}{[n]_{p,q}} p^k f \bigg(\frac {p^{n-k-1}[k]_{p,q}+ \alpha}{[n]_{p,q}+ \beta}\bigg)+\frac{[k]_{p,q}}{[n]_{p,q}} q^{n-k} f\bigg(\frac{p^{n-k}[k]_{p,q}+ \alpha}{[n]_{p,q}+ \beta}\bigg)- f \bigg (\frac {p^{n-k}[k]_{p,q}+ \alpha}{[n]_{p,q}+ \beta}\bigg).$\\

From \eqref{e14} it is clear that each $ \psi_k(x)$ is non-negative on $[0,1]$ for $0< q < p \leq 1 $ and, thus, it suffices to show that each $a_k$ is non-negative.\\

\noindent Since $f$ is convex on $[0,1]$ then for any $t_0 ,t_1$  and $\lambda \in [0,1]$ it follows that\\

$$ f(\lambda t_0 + (1-\lambda)t_1) \leq \lambda f(t_0) + (1-\lambda)f(t_1).$$\\
If we choose $t_0 = \frac {p^{n-k}[k]_{p,q}+ \alpha}{[n]_{p,q}+ \beta}$, $ t_1 = \frac {p^{n-k-1}[k]_{p,q}+ \alpha}{[n]_{p,q}+ \beta},$ and\\

$\lambda = \frac {[k]_{p,q}}{[n]_{p,q}}q^{n-k},$ then $t_0 ,t_1$ $\in $ $[0,1]$ and $\lambda \in (0,1)$ for $1 \leq k \leq n-1,$ and we deduce that \\
$$a_k = \lambda f( t_0) + (1-\lambda)f(t_1)- f(\lambda t_0 + (1-\lambda)t_1)\geq0 $$\\
Thus $ S_{n-1,p,q}(f;x)\geq S_{n,p,q}(f;x).$\\

We have equality for $x=0$ and $x=1,$ since the Bernstein polynomials interpolate $f$ on these end points.The inequality will be strict for $ 0 < x < 1 $ unless when $f$ is linear in each of the intervals between consecutive knots

\begin{equation*}
\frac{p^{n-k-1}[k]_{p,q} + \alpha}{[n]_{p,q} + \beta},~~       0 \leq k \leq n-1,
\end{equation*}

then we have\\
$$ S_{n-1,p,q}(f;x) =  S_{n,p,q}(f;x)$$       for $ 0 \leq x \leq 1.$\\
\section{A Global Approximation theorem}
In this section, we establish a global approximation theorem by means of Ditzian-Totik modulus of smoothness and Voronovskaja type approximation result.\\
In order to prove our next result, we recall the definitions of the Ditzian-Totik first order modulus of smoothness and the K-functional. Let $\phi(x) = \surd{x(1-x)}$ and $f \in C[0,1]$. The first order modulus of smoothness is given by\\
\begin{equation}\label{e15}
\omega_\phi(f;t) = \sup \limits _ {0< h \leq t} \Big\{\big|f(x + \frac {h\phi(x)}{2}) - f(x - \frac {h\phi(x)}{2})\bigg| , x \pm \frac {h\phi(x)}{2} \in [0,1]\Big\}\\
\end{equation}
The corresponding k-functional to \eqref{e15} is defined by\\

$k_\phi(f;t) =   \inf\limits_{g \in W_\phi [0,1]}\Big\{\|f - g\| + t \|\phi g^{\prime}\|\Big\} $    $(t > 0),$\\

where $W_\phi[0,1] = \{g: g \in AC_{loc}[0,1] , \|\phi g^{\prime}\|< \infty\} $ and $g \in AC_{loc}[0,1]$ means that $g$ is absolutely continuous on every interval $[a,b]\subset[0,1]$. It is well known \cite{dzk1} that there exists a constant $C >0$ such that\\
\begin{equation}\label{ab}
k_\phi(f;t) \leq Cw_\phi(f;t).\\
\end{equation}

\begin{theorem}
	Let $f \in C[0,1]$ and $\phi(x) = \surd{x(1-x)},$ then for every $x \in [0,1]$ we have \\
	$\bigg| S_{n,p,q}(f;x)- f(x)\bigg| \leq  C \omega_\phi\bigg(f;\frac{[n]_{p,q}}{\surd([n]_{p,q}+\beta)}\bigg)$
	where C is a constant independent of n and x.\\
\end{theorem}

\textbf{Proof:} Using the representation \\
$$ g(t) = g(x)+\int_{x}^{t} g^{\prime}(u) du ,$$\\

we get \\
\begin{equation}\label{e15}
\bigg| S_{n,p,q}(g;x)- g(x)\bigg| = \bigg| S_{n,p,q}\bigg(\int_{x}^{t} g^{\prime}(u) du ;x\bigg)\bigg|.\\
\end{equation}
For any $x \in (0,1)$ and $ t \in [0,1],$ we find that \\
\begin{equation}\label{e16}
\bigg |\int_{x}^{t} g^{\prime}(u) du\bigg |\leq \| \phi g^{\prime}\| \bigg|\int_{x}^{t} \frac {1}{\phi(u)}du \bigg|
\end{equation}
Further, \\
\begin{align}\label{e17}
\bigg| \int_{x}^{t} \frac {1}{\phi(u)}du\bigg| &= \bigg|\int_{x}^{t} \frac {1}{\surd u(1-u)} du\bigg| \notag\\
&\leq  \bigg| \int_{x}^{t} \bigg( \frac{1}{\surd u} + \frac{1}{\surd 1-u}\bigg)du\bigg |\notag\\
&\leq 2 ( |\surd{t} -\surd {x}|+ | \surd {1-t} - \surd {1-x}|)\notag\\
& = 2|t - x |\bigg (\frac {1}{\surd t + \surd x} + \frac{1}{\surd {1-t} + \surd {1-x}}\bigg)\notag\\
& < 2 |t - x | \bigg (\frac {1}{\surd x} +\frac {1}{\surd {1-x}}\bigg) \leq \frac { 2\surd 2 |t - x |}{\phi(x)}
\end{align}
From (\ref{e15}) - (\ref{e17}) and using the Cauchy - Schwarz inequality, we obtain \\
\begin{align*}
| S_{n,p,q}(g;x)- g(x)|&< 2\surd2 \|\phi g^{\prime}\| \phi^{-1}(x)S_{n,p,q}(|t - x|;x)\\
&\leq 2\surd2 \|\phi g^{\prime}\| \phi^{-1}(x)(S_{n,p,q}((t - x)^{2};x))^{\frac {1}{2}}.
\end{align*}

Using lemma (2.2), we get
\begin{align*}
| S_{n,p,q}(g;x)- g(x)| \leq \frac {2\surd2[n]_{p,q}}{\surd([n]_{p,q}+\beta} \|\phi g^{\prime}\|.
\end{align*}
Now using the above inequality we can write \\
\begin{align*}
| S_{n,p,q}(f;x)- f(x)|&\leq | S_{n,p,q}(f-g ;x)| + |f(x) - g(x)| + | S_{n,p,q}(g;x)- g(x)|\\
& \leq 2\surd 2 \bigg(\|f - g\|+ \frac {[n]_{p,q}}{\surd([n]_{p,q}+\beta)} \|\phi g^{\prime}\|\bigg).
\end{align*}
Taking the infimum on the right hand side of the above inequality over all $g \in W_\phi[0,1],$ we get\\
\begin{equation}\label{e18}
| S_{n,p,q}(f;x)- f(x)| \leq CK_\phi \bigg(f;\frac{[n]_{p,q}}{\surd([n]_{p,q}+\beta)} \bigg).
\end{equation}
Using equation $(\ref{ab})$ this theorem is proven.
%\begin{align*}
% \vert \int_{x}^{t}\vert(t-u)\vert \vertf^{,,}(u)_f^{,,}(x)\vert du \vert \leqslant 2 \Vert f^{,,} - g \Vert (t-x)^2 + 2\Vert \phi g^{\prime} \Vert \phi^{-1}(x)\vert t-x \vert^{3}\\
%\end{align*}
where $g\in W_\phi [0,1]$.On the other hand, for any  $ m = 1,2,.......$ and $ 0< q < p \leqslant 1$, there exists a constant $C_m > 0 $ such that\\
\begin{equation}\label{e20}
\vert S_{n,p,q}((t-x)_{p,q}^{m};x)\vert \leqslant C_m \frac{\phi^2(x)[n]_{p,q}}{([n]_{p,q}+ \beta)^{\lfloor \frac{m+1}{2}\rfloor}},
\end{equation}
where $x\in [0,1]$ and $\lfloor a \rfloor $ is the integral part of $a \geq 0.$\\

Throughout this proof, C denotes a constant not necessarily the same at each occurrence.\\

Now combining (\ref{e18}) -(\ref{e20}) and applying lemma (2.2), the cauchy-schwarz inequality,\\
We get \\
$\bigg| S_{n,p,q}(f;x)- f(x)\frac {p^{n-1}[n]_{p,q} - 2\alpha \beta}{2([n]_{p,q}+ \beta)^2}f^{\prime\prime}(x)\bigg|$\\
\begin{align*}
&\leq 2 \|f^{\prime\prime} - g \|S_{n,p,q}((t- x )^{2};x) + 2 \|\phi g^{\prime} \|\phi^{-1}(x)S_{n,p,q}(|t - x|^3;x)\\
&\leq2\|f^{\prime\prime} - g\|\frac{\phi^2(x)[n]_{p,q}}{([n]_{p,q}+ \beta)} + 2\|\phi g ^{\prime} \|\phi^{-1}(x)\{S_{n,p,q}((t- x )^{2};x)\}^{\frac{1}{2}}\{S_{n,p,q}((t- x )^{4};x)\}^{\frac{1}{2}}\\
&\leq2\|f^{\prime\prime} - g \|\frac{\phi^2(x)[n]_{p,q}}{([n]_{p,q}+ \beta)} + 2\frac{C}{([n]_{p,q}+ \beta)}\|\phi g ^{\prime} \|\frac{ \phi(x)[n]_{p,q}}{([n]_{p,q}+ \beta)^\frac{1}{2}}\\
&\leq \frac{C[n]_{p,q}}{([n]_{p,q}+ \beta)}\{\phi^2(x)\|f^{\prime\prime} - g \| + ([n]_{p,q}+ \beta)^\frac{-1}{2}\phi(x) \| \phi g ^{\prime}\|\}.\\
\end{align*}
since $\phi^2(x)\leq \phi(x)\leq 1,x\in[0,1],$ We obtain\\
\begin{align*}
\bigg|([n]_{p,q}+ \beta)^2[S_{n,p,q}(f;x)- f(x)]-\frac {p^{n-1}[n]_{p,q}-2\alpha\beta}{2}\phi^2(x)f^{\prime\prime}(x)\bigg| & \leq C\{\|f^{\prime\prime}-g\| \\
&+([n]_{p,q}+ \beta)^\frac{-1}{2}\phi(x)\|\phi g^{\prime}\|\}.\\
\end{align*}
Also, the following inequality can be obtained:\\
\subsection{Voronovskaja type theorem}
Using the first order Ditzian-Totik modulus of smoothnes, we prove a quantitative Voronovskaja type theorem for the $(p,q)$-Bernstein operators.
For any $f$ $\in C^2[0,1],$ the following inequalities holds:\\
\begin{equation}
\vert([n]_{p,q}+ \beta)[ S_{n,p,q}(f;x)- f(x)] - \frac{p^{n-1}- 2 \alpha\beta}{2} \phi^2(x)f^{\prime\prime}(x)\vert \leqslant C\omega_\phi(f^{\prime\prime}\phi(x)n^{\frac{-1}{2}}),
\end{equation}
\begin{equation}
\vert([n]_{p,q}+\beta)[ S_{n,p,q}(f;x)- f(x)] - \frac{p^{n-1}- 2 \alpha\beta}{2}\phi^2(x)f^{\prime\prime}(x)\vert \leqslant C\phi(x)\omega_\phi(f^{\prime\prime},n^{\frac{-1}{2}}),
\end{equation}
where C is a positive constant.\\

\textbf{Proof:} Let $f \in C^2[0,1]$ be given and $t,x \in[0,1]$ using Taylor's expansion, we have \\
\begin{equation}
f(t)-f(x) = (t-x)f^{\prime}(x)+ \int_{x}^{t}(t-u)f^{\prime\prime}(u)du
\end{equation}
Therefore
\begin{align*}
f(t) -f(x) - (t-x)f^{\prime}(x)- \frac{1}{2}(t-x)^2f^{\prime\prime}(x) &=  \int_{x}^{t}(t-u)f^{\prime\prime}(u)du - \int_{x}^{t}(t-u)f^{\prime\prime}(x)dx \\
&= \int_{x}^{t}(t-u)[f^{\prime\prime}(u) - f^{\prime\prime}(x)]du
\end{align*}
in view of lemma (2.2), we get\\
\begin{equation}
\bigg| S_{n,p,q}(f;x)- f(x) - \frac {p^{n-1}[n]_{p,q} - 2\alpha \beta}{2([n]_{p,q}+ \beta)^2}\phi^2(x)f^{\prime\prime}(x)\bigg| \leq  S_{n,p,q} \bigg(\bigg| \int_{x}^{t}|(t-u)||f^{\prime\prime}(u)- f^{\prime \prime}(x)|du \bigg|;x \bigg).
\end{equation}
The quantity $ |\int_{x}^{t} |f^{\prime\prime}(u)- f^{\prime\prime}(x)||(t-u)|du|$ was estimated in [ ],p- , as follows:\\
\begin{equation}
\bigg|\int_{x}^{t}f^{\prime\prime}(u)- f^{\prime \prime}(x)||t-u|du \bigg| \leq 2 \|f^{\prime\prime} - g \|(t - x)^2 + 2\| \phi g^{\prime}\| \phi^{-1}(x)|t - x|^3,
\end{equation}
%\begin{align*}
% \vert \int_{x}^{t}\vert(t-u)\vert \vertf^{,,}(u)_f^{,,}(x)\vert du \vert \leqslant 2 \Vert f^{,,} - g \Vert (t-x)^2 + 2\Vert \phi g^{\prime} \Vert \phi^{-1}(x)\vert t-x \vert^{3}\\
%\end{align*}
where $g\in W_\phi [0,1]$ . On the other hand, for any  $ m = 1,2,.......$ and $ 0< q < p \leqslant 1$, there exists a constant $C_m > 0 $ such that\\
\begin{equation}
\vert S_{n,p,q}((t-x)_{p,q}^{m};x)\vert \leqslant C_m \frac{\phi^2(x)[n]_{p,q}}{([n]_{p,q}+ \beta)^{\lfloor \frac{m+1}{2}\rfloor}}
\end{equation}

where $x\in [0,1]$ and $\lfloor a \rfloor $ is the integral part of $a \geq 0.$\\

Throughout this proof, C denotes a constant not necessarily the same at each occurrence.\\

Now combining (8.4) -(8.5) and applying lemma (2.2), the cauchy-schwarz inequality,\\
We get \\
$\bigg| S_{n,p,q}(f;x)- f(x)\frac {p^{n-1}[n]_{p,q} - 2\alpha \beta}{2([n]_{p,q}+ \beta)^2}f^{\prime\prime}(x)\bigg|$\\
\begin{align*}
&\leq 2 \|f^{\prime\prime} - g \|S_{n,p,q}((t- x )^{2};x) + 2 \|\phi g^{\prime} \|\phi^{-1}(x)S_{n,p,q}(|t - x|^3;x)\\
&\leq2\|f^{\prime\prime} - g\|\frac{\phi^2(x)[n]_{p,q}}{([n]_{p,q}+ \beta)} + 2\|\phi g ^{\prime} \|\phi^{-1}(x)\{S_{n,p,q}((t- x )^{2};x)\}^{\frac{1}{2}}\{S_{n,p,q}((t- x )^{4};x)\}^{\frac{1}{2}}\\
&\leq2\|f^{\prime\prime} - g \|\frac{\phi^2(x)[n]_{p,q}}{([n]_{p,q}+ \beta)} + 2\frac{C}{([n]_{p,q}+ \beta)}\|\phi g ^{\prime} \|\frac{ \phi(x)[n]_{p,q}}{([n]_{p,q}+ \beta)^\frac{1}{2}}\\
&\leq \frac{C[n]_{p,q}}{([n]_{p,q}+ \beta)}\{\phi^2(x)\|f^{\prime\prime} - g \| + ([n]_{p,q}+ \beta)^\frac{-1}{2}\phi(x) \| \phi g ^{\prime}\|\}\\
\end{align*}
since $\phi^2(x)\leq \phi(x)\leq 1,x\in[0,1],$ We obtain\\
\begin{equation}
\bigg|([n]_{p,q}+ \beta)^2[S_{n,p,q}(f;x)- f(x)] - \frac {p^{n-1}[n]_{p,q} - 2\alpha \beta}{2}\phi^2(x)f^{\prime\prime}(x)\bigg| \leq C\{\|f^{\prime\prime} - g\| + ([n]_{p,q}+ \beta)^\frac{-1}{2}\phi(x) \| \phi g ^{\prime}\|\}\\
\end{equation}
Also, the following inequality can be obtained:\\
\begin{equation}
\bigg|([n]_{p,q}+ \beta)^2[S_{n,p,q}(f;x)- f(x)] - \frac {p^{n-1}[n]_{p,q} - 2\alpha \beta}{2}\phi^2(x)f^{\prime\prime}(x)\bigg| \leq C\phi(x)\{\|f^{\prime\prime} - g\| + ([n]_{p,q}+ \beta)^\frac{-1}{2}\| \phi g ^{\prime}\|\}
\end{equation}
Taking the infimum on the right - hand side of the above relations over $g\in W_\phi[0,1],$ we get\\
\begin{equation}
\bigg|([n]_{p,q}+ \beta)^2[S_{n,p,q}(f;x)- f(x)] - \frac {p^{n-1}[n]_{p,q} - 2\alpha \beta}{2}\phi^2(x)f^{\prime\prime}(x)\bigg| \leq {C\phi(x)K_\phi( f^{\prime\prime};([n]_{p,q}+ \beta)^\frac{-1}{2})}{CK_\phi(f^{\prime\prime};\phi(x)([n]_{p,q}+ \beta)^\frac{-1}{2})},
\end{equation}
Using (8.9) and(7.2) the theorem is proved.\\
%\end{theorem}
\newpage

\section{Graphical Analysis}

With the help of Matlab, we show comparisons and some illustrative graphics for the
convergence of operators $(\ref{ee2})$ to the function $f(x)=1+x^3~ sin(14x)$  under different
parameters.\\

From figure \ref{f1}(a), it can be observed that as the value the $q ~\text{and}~ p $ approaches towards $1$ provided $0<q<p\leq1$, $(p, q)$-Bernstein Stancu operators given by $(\ref{ee2})$ converges towards the function.\\

From figure \ref{f1}(a) and (b), it can be observed that for $\alpha=\beta=0,$ as the value the $n$ increases, $(p, q)$-Bernstein Stancu operators given by \ref{ee2} converges towards the function $f(x)=1+x^3~ sin(14x)$.\\

Similarly from figure \ref{f2}(a), it can be observed that  for $\alpha=\beta=3,$ as the value the $q ~\text{and}~ p $ approaches towards $1$ provided $0<q<p\leq1$, $(p, q)$-Bernstein Stancu operators given by \ref{ee2} converges towards the function.\\

From figure \ref{f2}(a) and (b), it can be observed that as the value the $n$ increases, $(p, q)$-Bernstein Stancu operators given by $f(x)=1+x^3~ sin(14x)$ converges towards the function.

\begin{figure}[ht]\label{f1}
	\centering
	\subfigure[]
	{
		\includegraphics[height=4cm, width=6cm]{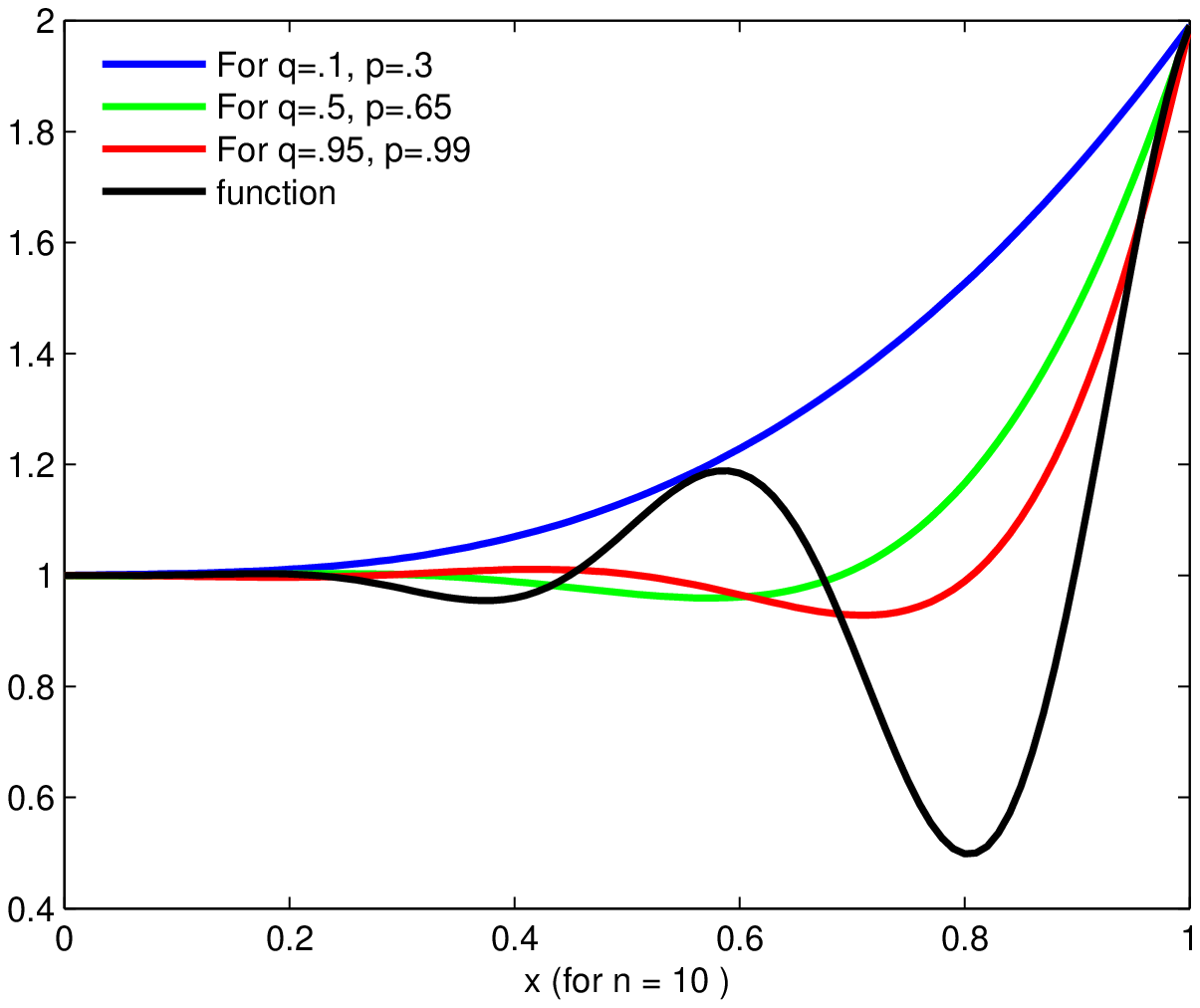}
		\label{fig:first_sub}
	}
	\subfigure[]
	{
		\includegraphics[height=4cm, width=6cm]{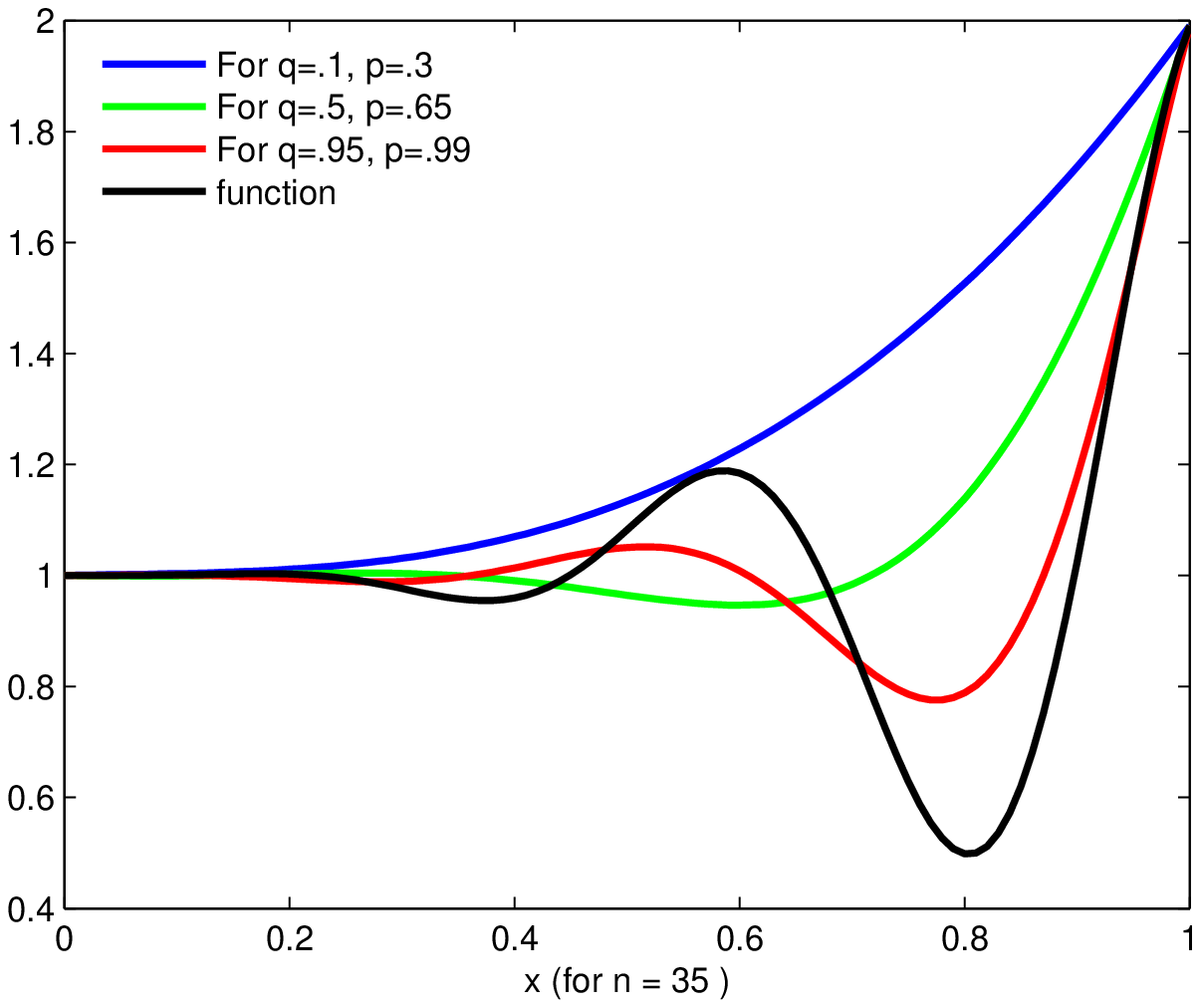}
		\label{fig:second_sub}
	}
	%\subfigure[Third caption]
	%    {
	%        \includegraphics[width=1.0in]{mkz_epi3.eps}
	%        \label{fig:third_sub}
	%    }
	\caption{$(p, q)$-Bernstein Stancu operators}\label{f1}
	%    \label{fig:sample_subfigures}
\end{figure}

\begin{figure}[ht]\label{f2}
	\centering
	\subfigure[]
	{
		\includegraphics[height=4cm, width=6cm]{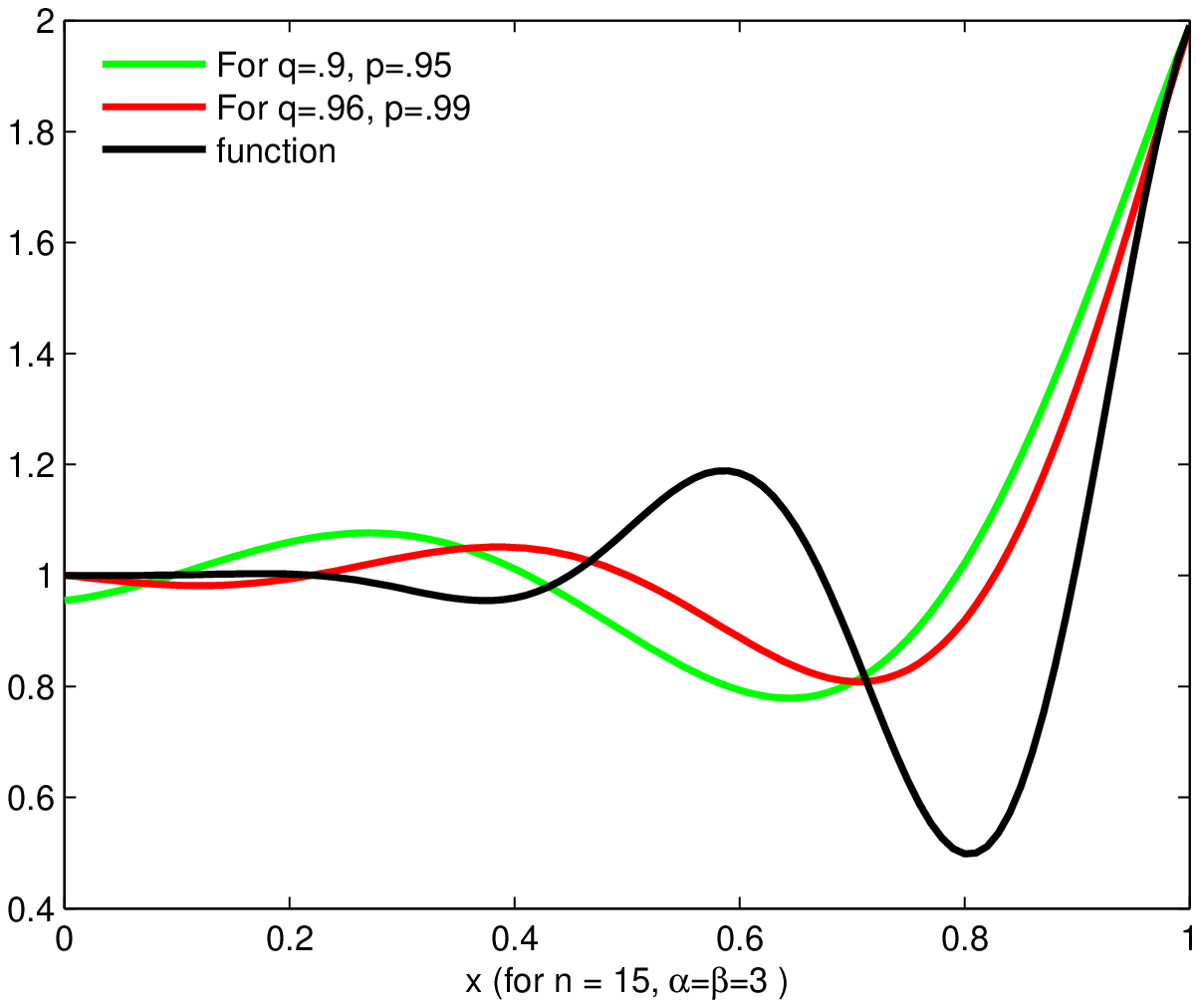}
		\label{fig:first_sub}
	}
	\subfigure[]
	{
		\includegraphics[height=4cm, width=6cm]{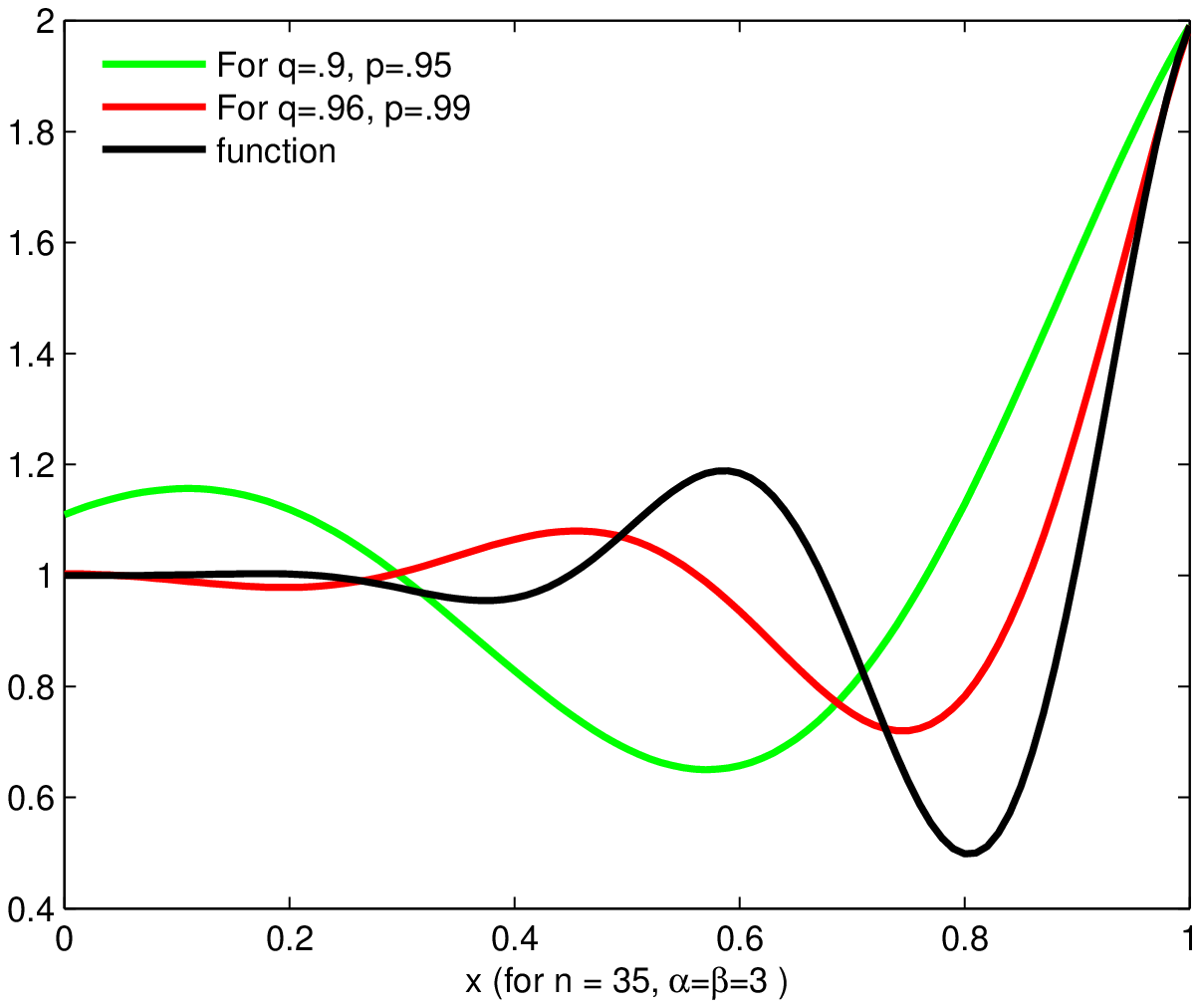}
		\label{fig:second_sub}
	}
	%\subfigure[Third caption]
	%    {
	%        \includegraphics[width=1.0in]{mkz_epi3.eps}
	%        \label{fig:third_sub}
	%    }
	\caption{$(p, q)$-Bernstein Stancu operators.}\label{f2}
	%    \label{fig:sample_subfigures}
\end{figure}

\end{document}